\documentclass{amsart}
\usepackage{latexsym,amsxtra,amscd,ifthen}
\usepackage{amsfonts}
\usepackage{verbatim}
\usepackage{amsmath}
\usepackage{amsthm}
\usepackage{amssymb}

\usepackage{enumerate}

\numberwithin{equation}{subsection}

\theoremstyle{plain}
\newtheorem*{theorem}{Theorem}
\newtheorem*{lemma}{Lemma}

\newtheorem*{corollary}{Corollary}

\theoremstyle{definition}
\newtheorem*{definition}{Definition}

\newtheorem*{remarks}{Remarks}



\DeclareMathOperator{\D}{{\sf D}} 
\DeclareMathOperator{\RHom}{RHom}

 \DeclareMathOperator{\Mod}{{\sf
Mod}}

\begin{document}

\title[Noetherian Hopf algebras]
{Noetherian Hopf algebras}

\author{Kenneth A. Brown}

\address{Department of Mathematics,
University of Glasgow, Glasgow G12 8QW, UK}

\email{kab@maths.gla.ac.uk}

\begin{abstract} This short survey article reviews our current state
of understanding of the structure of noetherian Hopf algebras. The
focus is on homological properties. A number of open problems are
listed.
\end{abstract}



\thanks{We thank Uli Kraehmer, Monica Macaulay, Catharina Stroppel
and James Zhang for helpful comments.}

\maketitle
\begin{center}{{\it To the memory of my teacher and friend Brian
Hartley}}\end{center}

\bigskip

\setcounter{section}{0}
\section{Introduction}
\label{intro}

\subsection{}
For the first 30 years after Hopf algebras were defined by H. Hopf
around 1940 the theory developed quite slowly. The publication of
Sweedler's monograph \cite{Sw} in 1969 quickened the pace, so that
understanding of the finite dimensional case in particular grew
considerably in the 1970s. But the tectonic plates really shifted
with the discovery of quantum groups \cite{D}, \cite{J} in the early
1980s, and the years since then have witnessed a massive expansion
in both the range of known examples and of our understanding of
them.

Many of these new examples of the past 25 years have been
noetherian algebras, so it makes sense to ask what features
noetherian Hopf algebras hold in common, and which aspects of the
finite dimensional theory extend to infinite dimensional
noetherian Hopf algebras. (We remark in passing that artinian Hopf
algebras give us nothing new, since every artinian Hopf algebra is
finite dimensional \cite{LiZ}.) Such an investigation was proposed
in the survey article \cite{Br1}, presented at an AMS meeting in
Seattle in 1997. The purpose of the present article is to review
what has happened since then: there have indeed been some
interesting and beautiful developments. As well as describing some
of these, I will list a number of questions which may help to
stimulate research on noetherian Hopf algebras over the next
decade.

\section{Definition and examples}
\label{defnex} \subsection{}\label{Hopfdef} All the algebras in
this paper will be defined over a field $k$ which for convenience
we shall always assume to be algebraically closed. To say that an
algebra $A$ is \emph{affine} means that $A$ is finitely generated
as an algebra. A Hopf algebra $H$ is an associative $k-$algebra
with a unit element, which is also equipped with
\begin{enumerate} \item a \emph{counit}; that is, an algebra
homomorphism $\varepsilon : H \longrightarrow k$; \item a
\emph{comultiplication}; that is, an algebra homomorphism $\Delta
: H \longrightarrow H \otimes_k H$, which we write using the
Sweedler notation: $\Delta(h) = \sum h_{1} \otimes h_2$ for $h \in
H$; \item an \emph{antipode}; that is, an algebra
antihomomorphism\footnote{An antihomomorphism is an algebra
homomorphism from $H$ to $H^{\mathrm{op}}$.} $S: H \longrightarrow
H.$
\end{enumerate}
This apparatus is required to satisfy a number of axioms
(essentially the duals of the axioms for a group). We won't list
these here as they can be found in all the standard references, for
example in \cite{Mo}, \cite{Sc}, \cite{Sw}. In addition, we'll
assume\footnote{By no means all results stated here require this
hypothesis, but we won't complicate matters by discussing details.}
throughout that
\begin{eqnarray} \textit{the antipode } S \textit{ is bijective.} \label{antipode1}
\end{eqnarray}
This hypothesis may in fact be vacuous - see (\ref{antip}) for a
discussion. We'll usually assume also that our Hopf algebras $H$
are \emph{left noetherian} - that is, all their left ideals are
finitely generated. Thanks to the antiautomorphism of $H$
gauranteed by (\ref{antipode1}), this is equivalent to $H$ being
right noetherian.

Recall that $H$ is said to be \emph{cocommutative} if
$\Delta^{\mathrm{op}}(h) := \sum h_2 \otimes h_1 = \Delta (h)$ for
all $h \in H.$ In the list of examples below we shall first review
the most important classes of cocommutative Hopf algebras (Exs. 1
and 2); then discuss the noetherian \emph{commutative} Hopf
algebras (Exs. 3) ; and then consider some classes of noetherian
Hopf algebras which may be neither cocommutative nor commutative
(Exs. 4-6).

\subsection{Examples} \label{examples} {\bf 1. Group algebras.} For
any group $G$, the group algebra $H = kG$ is a Hopf algebra, with
$\varepsilon (g)=1$, $\Delta (g) = g \otimes g$ and $S(g) =
g^{-1}$ for $g \in G.$ By a variant of Hilbert's basis theorem due
to Philip Hall \cite[Corollary 10.2.8]{Pa}, $kG$ is noetherian if
$kG$ is polycyclic-by-finite, (where this means that $G$ has a
finite series $1 = G_0 \subseteq G_1 \subseteq \ldots \subseteq
G_n = G$ of subgroups, with $G_i \lhd G_{i+1}$ and $G_{i+1}/G_i$
cyclic or finite for $i = 0, \ldots , n-1.)$ It's easy to see that
if $kT$ is any noetherian group algebra then $T$ satisfies the
ascending chain condition on subgroups, but more than 50 years
after Hall proved his theorem it's still not known if $T$ has to
be polycyclic-by-finite. So we ask:

\medskip
\noindent {\bf Question A:} Let $kG$ be a noetherian group
algebra. Is $G$ polycyclic-by-finite?
\medskip

\noindent {\bf 2. Enveloping algebras.} Let $\mathfrak{g}$ be a
$k-$Lie algebra. Then the enveloping algebra $H =
\mathcal{U}(\mathfrak{g})$ is a Hopf algebra with $\varepsilon (x) =
0,$ $\Delta (x) = x \otimes 1 + 1 \otimes x$ and $S (x) = -x$ for $x
\in \mathfrak{g}.$ By (the proof of) the
Poincar$\acute{\mathrm{e}}$-Birkhoff-Witt theorem,
$\mathcal{U}(\mathfrak{g})$ is a filtered algebra whose associated
graded algebra is a commutative polynomial algebra in
$\mathrm{dim}_k(\mathfrak{g})$ indeterminates. Thus, if
$\mathrm{dim}_k(\mathfrak{g}) < \infty ,$
$\mathcal{U}(\mathfrak{g})$ is a noetherian domain. It doesn't seem
to be known whether there are any other examples:

\medskip
\noindent {\bf Question B:} Suppose that
$\mathcal{U}(\mathfrak{g})$ is noetherian. Is
$\mathrm{dim}_k(\mathfrak{g}) < \infty$?

\medskip
Over characteristic 0 Examples 1 and 2 are not far from the
complete story for cocommutative Hopf algebras - by theorems of
Cartier, Gabriel and Kostant \cite[Corollary 5.6.4(3) and Theorem
5.6.5]{Mo}, if $k$ has charactersitic 0 and $H$ is \emph{any}
cocommutative Hopf $k-$algebra (not necessarily noetherian), then
$H$ is a skew group algebra over the enveloping algebra of the Lie
algebra of primitive elements of $H$.\footnote{This needs the
algebraic closure of $k$.} However other examples can occur in
positive characteristic - see \cite[pages 82-83]{Mo}.

\medskip
\noindent {\bf 3. Commutative Hopf algebras.} The category of
commutative affine Hopf $k-$algebras is equivalent to the category
of algebraic groups over $k$ \cite[Corollary 1.7]{Sc}: if $G$ is
such a group then its coordinate ring $\mathcal{O}(G)$ is a Hopf
algebra, with $\varepsilon (f) = f(1_G),$ $\Delta (f)$ the function
in $\mathcal{O}(G) \otimes \mathcal{O}(G) \cong \mathcal{O}(G \times
G)$ defined by $\Delta(f)((x,y)) := f(xy)$ for $x,y \in G,$ and
$S(f)(x) := f(x^{-1})$ for $x \in G.$ And by a theorem of Molnar
\cite{M}, a commutative Hopf algebra is affine if and only if it is
noetherian.

\medskip
In contrast to the above examples, quantum groups are neither
commutative nor cocommutative. Speaking crudely, these split into
two families, 4(i) and 4(ii) below, which are, respectively,
deformations of some of the algebras in Examples 2 and Examples 3.

\noindent {\bf 4. Quantum groups.} (i) {\bf Quantized enveloping
algebras} The key examples of these are deformations of the
enveloping algebras of semisimple Lie algebras. For each finite
dimensional semisimple Lie algebra $\mathfrak{g}$ and each
non-zero scalar $q \in k$, (avoiding a few ``bad'' values), $H =
\mathcal{U}_q(\mathfrak{g})$ is a noncommutative noncocommutative
noetherian Hopf $k-$algebra.

(ii) {\bf Quantized coordinate rings} For each semisimple
algebraic $k-$group $G$ and non-zero scalar $q$ (again avoiding a
few values), the quantized coordinate ring $H:=\mathcal{O}_q(G)$
is a deformation of the classical coordinate ring of $G.$ It is a
noncommutative noncocommutative noetherian Hopf algebra.

There are many references where details of the definitions and
basic properties of these algebras in Examples 4 can be found -
see, for example, \cite{Ja}, \cite{Jo}, \cite{BG2}. For $H$ in
either of the above classes, there is a fundamental dichotomy
determined by the value of the deformation parameter $q$: namely,
\begin{eqnarray} \label{PII}& H \textit{ is a finite module over its centre}\\ \nonumber & \textit{ if and only if } q \textit{ is a
root of 1 in } k. \end{eqnarray}

\medskip

\noindent {\bf 5. Hopf algebras satisfying a polynomial identity.}
For the definition of \emph{a ring satisfying a polynomial
 identity}, see for example \cite{MR}. The dichotomy (\ref{PII}) just identified for quantum groups can be
examined for the other example classes listed above. Thus a group
algebra $kG$ is a noetherian polynomial identity algebra if and only
if $G$ is a finitely generated abelian-by-finite group
\cite[Corollaries 5.3.8, 5.3.10]{Pa}. And the enveloping algebra
$\mathcal{U}(\mathfrak{g})$ of a finite dimensional Lie algebra
$\mathfrak{g}$ satisfies a polynomial identity if and only if
$\mathfrak{g}$ is abelian or $k$ has positive characteristic
\cite{La}, \cite{Z}. Prompted by this rather weak evidence, we ask
(i) below:

\medskip

\noindent {\bf Question C:} (i) Suppose that $H$ is a noetherian
Hopf algebra satisfying a polynomial identity. Is $H$ a finite
module over its centre?

(ii) (Wu, Zhang, \cite{WZ}) Let $H$ be as in (i). Is $H$ affine?

(iii) Let $H$ be an affine Hopf algebra satisfying a polynomial
identity. Is $H$ noetherian?

\medskip
Molnar's characterisation \cite{M} of commutative noetherian Hopf
algebras gives some support to (ii) and (iii). In \cite[Question
B]{Br1} I asked whether every affine noetherian PI Hopf algebra
was a finite module over a commutative normal sub-Hopf algebra.
(For the meaning of \emph{normal} here, see \cite[Definition
3.4.1]{Mo}.) It was noted by Gelaki and Letzter in \cite{GL} that
this is not the case, but their example does not rule out the
following refinement:

\medskip
\noindent {\bf Question D:} Suppose that $H$ is an affine
noetherian Hopf algebra satisfying a polynomial identity. Is $H$ a
finite module over a commutative normal right co-ideal subalgebra?

\medskip
This is true for all the PI algebras in the classes 1, 2 and 4.

We introduce the following class primarily so as to include factor
Hopf algebras of Examples 4(ii):

\noindent {\bf 6. Filtered algebras.} Let $H$ be a Hopf $k-$algebra.
We'll say that $H$ is \emph{normally $\mathbb{N}-$filtered} if $H =
\cup_{i \geq 0} H_i,$ with $H_0 = k$ and each $H_i$ a finite
dimensional $k-$vector space with $H_i H_j \subseteq H_{i+j}$ for
all $i,j$, such that the associated graded algebra $\mathrm{gr}(H)$
is connected graded noetherian, and so that every graded prime
factor ring of $\mathrm{gr}(H)$ is either $k$, or contains a
homogeneous normal element of positive degree.

\section{Motivation: finite dimensional Hopf algebras} \label{fdim}
In the subsequent sections we'll consider generalisations of the
classical facts about finite dimensional Hopf algebras which we
recall in (\ref{frob}) and (\ref{integral}).
\subsection{Frobenius algebras} \label{frob} Recall that a finite dimensional algebra $A$ is a \emph{Frobenius
algebra} if it admits a bilinear form $\phi : A \times A
\longrightarrow k$ which is non-degenerate (meaning that
$\phi(x,A) \neq 0 \neq \phi(A,x)$ for all $x \in A \setminus
\{0\}$), and associative (meaning that $\phi(xh,y) = \phi(x,hy)$
for all $x,y,h \in A$). Notice that this makes $A$ isomorphic to
its vector space dual $A^*$ as left and right $A-$module, so that
in particular $A$ is an injective left and right $A-$module - in
other words, $A$ is \emph{quasi-Frobenius}.

In 1969 Larson and Sweedler proved the following fundamental
theorem:

\begin{theorem} {\rm (Larson, Sweedler, \cite{LS})} Let $H$ be a finite dimensional Hopf algebra. Then
$H$ is a Frobenius algebra.
\end{theorem}

\subsection{Integrals} \label{integral} The left $H-$module
isomorphism of $H$ and $H^*$ implies that $H$ contains a unique
ideal $\int^l_H$ with $\mathrm{dim}_k(\int^l_H) = 1$ and $hx =
\varepsilon (h)x$ for all $x \in \int^l_H$. The ideal $\int^l_H$ is
called the \emph{left integral} of $H$. In a similar way $H$ has a
\emph{right integral} $\int^r_H.$ If $\int^l_H = \int^r_H$, $H$ is
called \emph{unimodular}. For example if $G$ is a finite group then
$H = kG$ is unimodular with $\int^l_H = \int^r_H = \sum_{g \in G}g.$

\section{Injective dimension} \label{injdim}
\subsection{}Self-injective algebras are artinian \cite[Proposition XIV.3.1]{St}, so
it's clear that Theorem \ref{frob} doesn't generalise directly to
infinite dimensional algebras. On the other hand, it's easy to see
that, when $k$ has characteristic 0, commutative affine Hopf
$k-$algebras are \emph{regular} - that is, they have finite global
(projective) dimension. (This is essentially because in
characteristic 0 commutative Hopf algebras are semiprime by a
theorem of Cartier \cite[Theorem 11.4]{Wa}, and the regular action
of the group $G$ defines automorphisms mapping any given maximal
ideal of $\mathcal{O}(G)$ to any other.) More generally, over any
field, commutative noetherian Hopf algebras are \emph{Gorenstein} -
that is, they have finite injective dimension, \cite[Proposition
2.3]{Br1}. Now \emph{any} commutative affine Gorenstein (or, \emph{a
fortiori}, regular) algebra has injective dimension equal to the
``size'' $d$ of the algebra, \cite[Theorem 21.8]{Eis}. Here,
``size'' means the Gelfand-Kirillov dimension, $\mathrm{GKdim}(-)$,
or Krull dimension (which are equal for a commutative affine
algebra). The \emph{Gelfand-Kirillov dimension} of an affine algebra
$A$ is a measure of its rate of growth; it has many attractive
properties, \cite{KL}, but unfortunately is often infinite.
\emph{Krull dimension}, on the other hand, is \emph{always} defined
for a noetherian algebra, but its use often involves difficult
technical problems. In any case, it seems that the correct way to
impose the relevant ``size'' constraints in a noncommutative setting
may be to demand more stringent homological conditions than simply
having finite injective dimension. The relevant definitions are
introduced in the next paragraph.

\subsection{Homological definitions} \label{hom} Useful sources for
the basic facts concerning the following ideas are \cite{ASZ1},
\cite{BG1}, \cite{Lev}. A simple but key point to appreciate when
considering (b) and (f) is that, for say a \emph{left} $A-$module
$M$, $\mathrm{Ext}_A^i(M,A)$ is a \emph{right} $A-$module via the
right action on $A$.
\begin{definition} Let $A$ be a ring.
\begin{enumerate}\item The \emph{grade} of the left $A-$module $M$ is
$$ j(M) := \mathrm{inf}\{j : \mathrm{Ext}_A^j(M,A) \neq 0 \}. $$
\item $A$ satisfies the \emph{Auslander condition} if, for every
noetherian left or right $A-$module $M$ and for all $i \geq 0,$
$j(N) \geq i$ for all submodules $N$ of $\mathrm{Ext}_A^i(M,A).$
\item The ring $A$ is \emph{Auslander-Gorenstein} if it is
noetherian, satisfies the Auslander condition, and has finite
right and left injective dimensions (which are then equal by a
theorem of Zaks \cite{Za}). \item If $A$ is Auslander-Gorenstein
and has finite global dimension then it is called
\emph{Auslander-regular}. \item The ring $A$ is
\emph{Cohen-Macaulay} (with respect to $\mathrm{GKdim}$) if, for
all non-zero noetherian $A-$modules $M$,
$$ j(M) + \mathrm{GKdim}(M) = \mathrm{GKdim}(A). $$
\item Suppose that $A$ is a noetherian Hopf $k-$algebra. Then $A$ is
$ AS-$\emph{Gorenstein} if it has right and left injective dimension
$d < \infty$, and $\mathrm{Ext}_A^i (k,A) = \delta_{id}k,$ where the
module $k$ is as usual the trivial (right or left) $A-$module, with
$A$ acting through $\varepsilon$.
\item The Hopf algebra $A$ is $AS-$\emph{regular} if it is $AS-$Gorenstein
and has finite global dimension.
\end{enumerate}
\end{definition}

These definitions are closely connected, at least for noetherian
Hopf algebras:

\begin{lemma}(\cite[Lemma 6.1]{BZ}) Let $H$ be a noetherian Hopf $k-$algebra. If $H$ is
Auslander-Gorenstein and Cohen-Macaulay, then $H$ is
$AS-$Gorenstein.
\end{lemma}

\subsection{Injective dimension of Hopf algebras}\label{injdim} In \cite[3.1]{Br1} and also in
\cite[1.15]{BG1} we asked whether every noetherian Hopf algebra has
finite injective dimension. This question remains open, so we
restate it here, taking the opportunity to refine it in the light of
evidence gathered in the last decade:
\medskip

\noindent {\bf Question E:} Is every noetherian Hopf algebra
$AS-$Gorenstein?

\medskip
At the time of writing, the answer is ``yes'' for all known
noetherian Hopf algebras. In particular, the algebras listed in
(\ref{examples}) are all $AS-$Gorenstein. Detailed proofs for
classes 1, 2 and 4(i) can be found in \cite[$\S$6]{BZ}; see
\cite{GZ} for class 4(ii). The proof for class (\ref{examples})6
given in \cite{LWZ2} is different in flavour; we discuss it briefly
in Remark (\ref{dual2})(b). The most striking of these positive
cases for Question E is class (\ref{examples})5, affine noetherian
Hopf algebras satisfying a polynomial identity - the result is a
theorem of Wu and Zhang which is both beautiful and technical. In
fact, at least formally, they prove a bit more:

\begin{theorem} {\rm (Wu, Zhang \cite{WZ})} Every affine noetherian Hopf algebra satisfying a
polynomial identity is Auslander-Gorenstein and Cohen-Macaulay.
\end{theorem}

To illustrate the power of these homological properties we state a
sample non-homological corollary, the second (much deeper) part of
which follows from the theorem together with results of Stafford and
Zhang \cite{SZ}:

\begin{corollary} Let $H$ be as in the theorem.
\begin{enumerate}
\item \cite[Theorem 0.2(2)]{WZ} $H$ has a quasi-Frobenius
(artinian) ring of fractions. \item Suppose that $H$ has finite
global dimension. Then $H$ is a finite direct sum of prime rings,
and is a finite module over its centre.
\end{enumerate}
\end{corollary}

\subsection{Integrals of Hopf algebras}\label{integrals} While it
is perhaps not so surprising that finiteness of the injective
dimension should generalise from artinian to noetherian Hopf
algebras, it was very surprising - to me at least - that the idea of
the integral should do so also. Let $_{\varepsilon}k$ denote the
trivial left $H-$module. The key points are first, to think of
$\int^l_H$ in the artinian case as $\mathrm{Hom}_H(_{\varepsilon}k,
H)$; second, to regard $\mathrm{Hom}_H(_{\varepsilon}k, H)$ as the
case $i=0$ of $\mathrm{Ext}^i_H(_{\varepsilon}k, H)$; and third, to
recall that these $\mathrm{Ext}-$groups are $H-$bimodules, with
\emph{left} $H-$action induced by the right action on
$_{\varepsilon}k$ (and so trivial), and \emph{right} action induced
from the right action on $H.$ The definition is due to Lu, Wu and
Zhang:

\begin{definition} \cite[Definition 1.1]{LWZ} Let $H$ be an
$AS-$Gorenstein Hopf algebra of injective dimension $d$.
\begin{enumerate} \item The one-dimensional $k-$vector space and
$H-$bimodule $\mathrm{Ext}^d_H(_{\varepsilon}k,\, _| H)$ is called
the \emph{left integral} of $H$, denoted $\int^l_H.$ \item The
one-dimensional $k-$vector space and $H-$bimodule
$\mathrm{Ext}^d_H(k_{\varepsilon},H_|)$ is called the \emph{right
integral} of $H$, denoted $\int^r_H.$ \item $H$ is
\emph{unimodular} if $\int^l_H$ is right trivial as well as left
trivial.
\end{enumerate}
\end{definition}

One can show quite easily \cite[Lemma 1.3]{LWZ} that $H$ is
unimodular if and only if $\int^r_H$ is left trivial.

\subsection{The Nakayama automorphism} \label{Nak} As we saw in
(\ref{frob}), if $A$ is any Frobenius algebra (for example a
finite dimensional Hopf algebra) then $A^*$ is isomorphic to $A$
as left and as right $A-$module. But in general this is \emph{not}
an isomorphism of bimodules: in fact the correction is provided by
twisting the module on one side by a suitable algebra automorphism
$$ A^* \, \cong \; ^{\nu}A^1,$$
\cite[Theorem 2.4.1]{Ya}. Here, $^{\nu}A^1$ is the $A-A-$bimodule
which is left and right isomorphic to $A$, with $a.b.c := \nu
(a)bc$ for all $a,c \in A,$ for all $b \in \, ^{\nu}A^1$. In the
theory of Frobenius algebras, $\nu$ is called the \emph{Nakayama
automorphism} of $A$, well-defined up to an inner automorphism of
$A$. For many purposes - for instance, in representation theory -
it's useful to know $\nu$ explicitly. When $\nu = \mathrm{Id}$,
$A$ is called a \emph{symmetric algebra}.

Recall that if $H$ is \emph{any} Hopf algebra (not necessarily
finite dimensional) and $\pi : H \longrightarrow k$ is an algebra
epimorphism, the \emph{left winding automorphism} $\tau^l_{\pi}$
is the algebra automorphism
$$ \tau^l_{\pi} : H \longrightarrow H : h \mapsto \sum
\pi(h_1)h_2. $$ The \emph{right winding automorphism}
$\tau^r_{\pi}$ is defined by $\tau^r_{\pi}(h) = \sum h_1\pi(h_2)$
for $h \in H.$ The Nakayama automorphism of a finite dimensional
Hopf algebra has the following description:

\begin{theorem} {\rm (Schneider, \cite[Proposition 3.6]{Sc})} Let $H$ be
a finite dimensional Hopf algebra and let $\pi : H \longrightarrow
k$ be the algebra epimorphism whose kernel is the right
annihilator of $\int^l_H$. Then the Nakayama automorphism $\nu$ of
$H$ is $S^2 \circ \tau^l_{\pi}.$
\end{theorem}

\section{Dualizing complexes} \label{dual} \subsection{} \label{dualdef} Theorem \ref{Nak}
generalises in a natural way to AS-Gorenstein Hopf algebras,
provided we work in the derived category, in particular using
concepts developed by Yekutieli \cite{Ye} and Van den Bergh
\cite{VdB2}. Recall that if $A$ is a noetherian algebra, a bounded
complex $R_A$ of $A-A-$bimodules (viewed as an object of the bounded
derived category $\D (A^e$-$\Mod)$ of $A-A-$bimodules) is a
\emph{rigid dualizing complex} over $A$ if
\begin{enumerate}
\item $R$ has finite injective dimension over $A$ and over $A^{\sf
op}$ respectively. \item $R$ is homologically finite over $A$ and
over $A^{\sf op}$ respectively. \item The canonical morphisms
$A\to \RHom_A(R,R)$ and $A\to \RHom_{A^{\sf op}}(R,R)$ are
isomorphisms in $\D (A^e$-$\Mod)$. \item  A dualising complex $R$
over $A$ is called {\it rigid} if there is an isomorphism
$$R \cong \RHom_{A^e}(A, R \otimes R^{\sf op})$$
in $\D (A^e$-$\Mod)$. (Here the $A-A-$bimodule structure of
$R\otimes R^{\sf op}$ comes from the left $A$-module structure of
$R$ and the left $A^{\sf op}$-module structure of $R^{\sf op}$).
\end{enumerate}

When such a complex $R$ exists it is unique, and
$\mathrm{RHom}_A(-, R)$ defines a duality - that is, a
contravariant equivalence - between the bounded derived categories
of left and right $A-$modules. For example, if $A$ is any finite
dimensional algebra then $R_A$ exists and is $A^*$. So if $A$ is a
Frobenius algebra,
$$ R_A = A^* \cong \; ^{\nu}A^1. $$

\subsection{}\label{dual2} If $M$ is an $A-$module and $d \in \mathbb{Z}$, we write
$M[d]$ for the complex which has $M$ moved $d$ places to the left
(from the $0$th place) and 0 elsewhere. We can now state a result
generalising this left-right duality from the finite-dimensional
case to noetherian AS-Gorenstein Hopf algebras:

\begin{theorem} \cite[Proposition 4.5]{BZ} Let $H$ be an AS-Gorenstein Hopf
algebra of injective dimension $d$. \begin{enumerate} \item $H$ has
rigid dualizing complex $^{\nu}H^1[d]$, for a certain algebra
automorphism $\nu$ of $H$.
\item The automorphism $\nu$, which we call the Nakayama
automorphism of $H$, is $S^2 \circ \tau^l_{\pi}$, where $\pi$ is the
epimorphism from $H$ to $H/(\mathrm{r-ann}(\int^l_H)).$
\end{enumerate}
\end{theorem}
\medskip

Naturally, we should ask the following question, which is probably
closely related to Question E:

\noindent {\bf Question F:} Does every noetherian Hopf algebra
have a rigid dualizing complex?

\begin{remarks} (a) It follows from the above that the Nakayama
automorphism and the integrals are crucial to the two-sided
structure of AS-Gorenstein Hopf algebras. The calculation of these
entities for classes (\ref{examples})1, 2 and 4 is not difficult and
has been carried out in \cite[$\S$6]{BZ}.

(b) The treatment \cite{LWZ2} of the normally
$\mathbb{N}-$filtered Hopf algebras of (\ref{examples})6 is the
reverse of that given here. Namely, one shows first that such an
$H$ \emph{has} a rigid dualizing complex satisfying a rather
natural additional property, and then deduces \emph{from this}
that $H$ is AS-Gorenstein. As this indicates, it seems that
Questions E and F are closely related.
\end{remarks}

\section{Applications of the dualizing complex} \label{apps}
\subsection{Poincar$\acute{\mathbf{e}}$ duality} \label{Poincare} For the
definition of the Hochschild homology groups $H_i(A,M)$ and
cohomology groups $H^i(A,M)$ of an $A-$bimodule $M$ we refer to
\cite[Chapter 9]{We}. Although classical
Poincar$\acute{\mathrm{e}}$ duality fails for noncommutative
noetherian Hopf algebras, it seems that it may be valid if we
allow twisting by the Nakayama automorphism. Combining Theorem
\ref{dual2} with a result of Van den Bergh \cite{VdB1} we obtain

\begin{theorem} Let $H$ be a noetherian AS-regular Hopf algebra of
global dimension $d$ with Nakayama automorphism $\nu$.
\begin{enumerate}
\item For every $A-$bimodule $M$ and all $i = 0, \ldots , d$
$$ H^i(H,\,M) \cong H_{d-i}(H,\, ^1M^{\nu}). $$
\item In particular,
$$ H^d(H,\, ^{\nu}H^1) \cong H/[H,H] \neq 0, $$
and
$$ H_d(H,\,^1H^{\nu}) \cong Z(A) \neq 0. $$
\end{enumerate}
\end{theorem}

\subsection{The antipode} \label{antipode} If $H$ is a Hopf
algebra (with bijective antipode as usual) then so is $H' :=
(H,\Delta^{op}, S^{-1}, \varepsilon )$, where $\Delta^{op} (h) :=
\sum h_2 \otimes h_1$ \cite[Lemma 1.5.11]{Mo}. If $H$ is noetherian
and AS-Gorenstein we can apply Theorem \ref{dual2}(b) to it
\emph{and} to $H'$. The latter case yields Nakayama automorphism
$$ \nu' = \tau^r_{\pi} \circ S^{-2}, $$
where $\tau^r_{\pi}$ is the \emph{right} winding automorphism
associated to the epimorphism $\pi : H \longrightarrow
H/\mathrm{r-ann}(\int^l_H)$; see (\ref{Nak}). However, the Nakayama
automorphism of $H$ is unique up to an inner automorphism, by the
uniqueness property of rigid dualizing complexes. Since the
underlying algebra for $H'$ is the same as for $H$, $\nu$ and $\nu'$
differ only by an inner automorphism, proving:

\begin{theorem} \cite[Corollary 4.6]{BZ} Let $H$ be a noetherian AS-Gorenstein Hopf
algebra. Then there exists an inner automorphism $\gamma$ such
that
$$ S^4 = \gamma \circ \tau^r_{\pi} \circ (\tau_{\pi}^l)^{-1}, $$
where $\tau^l_{\pi}$ and $\tau^r_{\pi}$ are the left and right
winding automorphisms given by the left integral of $H$.
\end{theorem}
\medskip

Of course we immediately ask:

\noindent {\bf Question G:} What is the inner automorphism
$\gamma$?
\medskip

When $H$ is finite dimensional $\gamma$ is conjugation by the
group-like element which is the character of the right structure
on $\int^l_{H^*}$, by a 1976 paper of Radford \cite{Ra}. This
suggests that the Hopf dual $H^{\circ}$ may feature in the answer
to Question G.

It's not hard to see that the maps $\gamma, \tau^l_{\pi}$ and
$\tau^r_{\pi}$ in the theorem commute with each other
\cite[Proposition 4.6]{BZ}. Moreover, when $H$ is a finite module
over its centre, $\tau^l_{\pi}$ and $\tau^r_{\pi}$ have finite
order \cite[Theorem 2.3(b)]{BZ2}. It follows that:

\begin{corollary} If $H$ is a noetherian Hopf algebra which is a
finite module over its centre, then some power of the antipode of
$H$ is inner.
\end{corollary}

\medskip
\noindent {\bf Question I:} Is the corollary true for an affine
noetherian Hopf algebra satisfying a polynomial identity?
\medskip

The antipode of a finite dimensional Hopf algebra has finite order
\cite{Ra}, and $S^2 = \mathrm{Id}$ for a commutative Hopf algebra,
\cite[Corollary 1.5.12]{Mo}, so it's natural to ask:

\noindent {\bf Question H:} If $H$ is as in the corollary, does
$S$ have finite order?

\section{Further questions}
\subsection{Finite global dimension} \label{global} The
possibility that \emph{all} noetherian Hopf algebras have finite
injective dimension, together with the motivating commutative,
cocommutative and finite dimensional cases, combine to suggest
that there may be natural structural conditions on a noetherian
Hopf algebra sufficient to guarantee other homological properties
such as finite global dimension. If we examine our favourite
classes of examples, at least three structural "indicators" of
infinite global dimension for a noetherian Hopf algebra $H$ become
apparent:

\begin{enumerate} \item $H$ is not semiprime; \item $H$ has a
finite dimensional Hopf subalgebra which is not semisimple; \item
$H$ has a finite dimensional irreducible module of dimension
divisible by the characteristic of $k$.
\end{enumerate}

Of course, more than one of these features can occur in the same
example; and easy group algebra examples show that (c) can happen
in a regular Hopf algebra. Nevertheless, in part to stimulate the
creation of more esoteric examples, we ask:
\medskip

\noindent {\bf Question J:} Suppose a noetherian Hopf algebra $H$
is not regular. Must at least one of (a), (b), (c) occur?
\medskip

If this seems too optimistic or too difficult, one might try:
\medskip

\noindent {\bf Question K:} Let $H$ be a noetherian domain and a
Hopf $k-$algebra, and suppose that $k$ has characteristic 0. Is
$H$ regular?
\medskip

Some (slight) positive evidence for Question J is given by

\begin{theorem} {\rm (Wu, Zhang \cite{WZ2})} Let $H$ be a noetherian
affine Hopf algebra satisfying a polynomial identity. Suppose that
$H$ is involutary - that is, $S^2 = \mathrm{Id}$. If neither (a)
nor (c) occurs for $H$, then $H$ is regular.
\end{theorem}

\subsection{Bijectivity of the antipode} \label{antip} Recall that
we've assumed throughout that our Hopf algebras have a bijective
antipode (\ref{antipode1}). Examples of Takeuchi \cite{Tak} show
that this hypothesis fails in general. However no example is known
of a \emph{noetherian} Hopf algebra whose antipode is \emph{not}
bijective, and we have the following theorem and final question:

\begin{theorem} {\rm (Skryabin, \cite{Sk})} If $H$ is a noetherian Hopf
algebra which is either semiprime or affine with a polynomial
identity, then its antipode is bijective.
\end{theorem}
\medskip

\noindent {\bf Question L:} (Skryabin) Let $H$ be a noetherian
Hopf algebra. Is the antipode $S$ bijective?


\begin{thebibliography}{10}


\bibitem[ASZ1]{ASZ1}
K. Ajitabh, S.P. Smith and J.J. Zhang, Auslander-Gorenstein rings,
Comm. Algebra {\bf 26} (1998), no. 7, 2159--2180.


\bibitem[Br1]{Br1}
K. A. Brown, Representation theory of Noetherian Hopf algebras
satisfying a polynomial identity, Trends in the representation
theory of finite-dimensional algebras (Seattle, WA, 1997), 49--79,
Contemp. Math., {\bf 229}, AMS, Providence, RI, 1998.

\bibitem[BG1]{BG1}
K. A. Brown and K. R. Goodearl, Homological aspects of Noetherian
PI Hopf algebras and irreducible modules of maximal dimension, J.
Algebra {\bf 198} (1997), 240--265.

\bibitem[BG2]{BG2}
K. A. Brown and K. R. Goodearl, \emph{Lectures on Algebraic Quantum
Groups}, Birkh{\"a}user, 2002.

\bibitem[BZ]{BZ}
K. A. Brown and J.J. Zhang, Dualizing complexes and twisted
Hochschild (co)homology for noetherian Hopf algebras, to appear in
J. Algebra; arXiv math/0603732.

\bibitem[BZ2]{BZ2}
K. A. Brown and J.J. Zhang, Prime regular Hopf algebras of
GK-dimension one, in preparation.


\bibitem[D]{D}
V.G. Drinfeld,
 Hopf algebras and the quantum Yang-Baxter equation, Sov. Math.
 Dokl. {\bf 32} (1985), 254-258.

\bibitem[Eis]{Eis}
D.Eisenbud, \emph{Commutative Algebra with a View toward Algebraic
Geometry}, Springer, 1995.

\bibitem[GL]{GL}
S. Gelaki and E. Letzter, An affine PI Hopf algebra not finite over
a normal commutative Hopf subalgebra, Proc. Amer. Math. Soc. {\bf
131} (2003), 2673--2679.

\bibitem[GZ]{GZ}
K. R. Goodearl and J. J. Zhang, Homological properties of quantized
coordinate rings of semisimple groups, Proc. London Math. Soc. {\bf
94} (2007), 647-671.



\bibitem[Ja]{Ja}
J.-C. Jantzen, \emph{Lectures on Quantum groups}, Graduate Studies
in Mathematics {\bf 6}, Amer. Math. Soc. 1996.

\bibitem[J]{J}
M. Jimbo, A $q-$difference analogue of $U(\mathfrak{g})$ and the
Yang-Baxter equation, Lett. Math. Physics, {\bf 10} (1985), 63-69.

\bibitem[Jo]{Jo}
A. Joseph, \emph{Quantum Groups and their Primitive Ideals},
Springer, 1995.


\bibitem[KL]{KL}
G. Krause and T.H. Lenagan, \emph{Growth of Algebras and
Gelfand-Kirillov Dimension (Revised Edition)}, Graduate Studies in
Mathematics {\bf 22}, Amer. Math. Soc. 2000.


\bibitem[LS]{LS}
R.G. Larson and M. Sweedler, An associative orthogonal bilinear form
for Hopf algebras, Amer. J. Math. {\bf 91} (1969) 75--94.

\bibitem[La]{La}
V.N. Latysev, Two remarks on PI-algebras, (Russian) Sibirsk. Mat. Z.
{\bf 4} 1963, 1120–1121.

\bibitem[Lev]{Lev}
T. Levasseur, Some properties of noncommutative regular graded
rings, Glasgow Math. J. 34 (1992), 277-300.


\bibitem[LiZ]{LiZ}
C.-H. Liu and J.J. Zhang, Artinian Hopf algebras are finite
dimensional, Proc. Amer. Math. Soc. 135 (2007), 1679-1680.

\bibitem[LWZ2]{LWZ2}
D.-M. Lu, Q.-S. Wu and J.J. Zhang, Hopf algebras with rigid
dualizing complex, Israel J. Math., to appear.



\bibitem[LWZ]{LWZ}
D.-M. Lu, Q.-S. Wu and J.J. Zhang, Homological integral of Hopf
algebras, Trans. Amer. Math. Soc. {\bf 359} (2007), 4945-4975.


\bibitem[MR]{MR}
J. C. McConnell and J. C . Robson, \emph{Noncommutative Noetherian
Rings}, Wiley, Chichester, 1987.

\bibitem[M]{M}
R.K. Molnar, A commutative Noetherian Hopf algebra over a field is
finitely generated, Proc. Amer. Math. Soc. {\bf 51} (1975), 501-502.

\bibitem[Mo]{Mo}
S. Montgomery, \emph{Hopf Algebras and their Actions on Rings},
         CBMS Regional Conference Series in
Mathematics, {\bf 82},
         Providence, RI, 1993.


\bibitem[Pa]{Pa}
D.S. Passman, \emph{The Algebraic Structure of Group Rings},
reprint of the 1977 original, Robert E. Krieger Publishing Co.,
Inc., Melbourne, FL, 1985.


\bibitem[Ra]{Ra}
D.E. Radford, The order of the antipode of a finite dimensional
Hopf algebra is finite, Amer. J. Math. (98) 1976, 333-355.


\bibitem[Sc]{Sc}
H.-J. Schneider, Lectures on Hopf Algebras, (Notes by Sonia
Natale), Trabajos de Matem{\'a}tica, 31/95, Facultad de
Matem{\'a}tica, Astronom{\'a}a y Física, C{\'o}rdoba, 1995;
www.famaf.unc.edu.ar/series/pdf/pdfBMat/BMat31.pdf

\bibitem[Sk]{Sk}
Skryabin, Serge, New results on the bijectivity of antipode of a
Hopf algebra, J. Algebra {\bf 306} (2006), 622--633.


\bibitem[St]{St}
B. Stenstrom, \emph{Rings of Quotients}, Springer-Verlag,
Heidelberg, 1975.


\bibitem[SZ]{SZ}
J. T. Stafford and J. J. Zhang, Homological properties of (graded)
Noetherian PI rings, J. Algebra {\bf 168} (1994), no. 3,
988--1026.

\bibitem[Sw]{Sw}
M. E. Sweedler, \emph{Hopf Algebras}, Benjamin, 1969.

\bibitem[Tak]{Tak}
M. Takeuchi, There exists a Hopf algebra whose antipode is not
injective, Sci. Papers College Gen. Ed. Univ. Tokyo, {\bf 21}
(1971), 127-130.

\bibitem[VdB1]{VdB1}
M. Van den Bergh, A relation between Hochschild homology and
cohomology for Gorenstein rings, Proc. Amer. Math. Soc. {\bf 126}
(1998), no. 5, 1345--1348. And Erratum to: "A relation between
Hochschild homology and cohomology for Gorenstein rings" Proc.
Amer. Math. Soc. {\bf 130} (2002), no. 9, 2809--2810.

\bibitem[VdB2]{VdB2}
M. Van den Bergh, Existence theorems for dualizing complexes over
noncommutative filtered and graded rings, J. Algebra {\bf 195}
(1997), 662-679.

\bibitem[Wa]{Wa}
W.C. Waterhouse, \emph{Introduction to Affine Group Schemes},
Graduate Texts in Mathematics 66, Springer, Berlin, 1979.


\bibitem[We]{We}
C. A. Weibel, \emph{An Introduction to Homological Algebra},
Cambridge Studies in Advanced Mathematics 38, Cambridge University
Press, Cambridge, 1994.

\bibitem[WZ]{WZ}
Q.-S. Wu and J.J. Zhang, Noetherian PI Hopf algebras are Gorenstein,
Trans. Amer. Math. Soc. {\bf 355} (2002), 1043-1066.

\bibitem[WZ2]{WZ2}
Q.-S. Wu and J.J. Zhang, Regularity of involutory PI Hopf algebras,
J. Algebra {\bf 256} (2002), 599--610.

\bibitem[Ya]{Ya}
K. Yamagata, Frobenius algebras, \emph{Handbook of Algebra}, Vol.1,
841-887, North-Holland, Amsterdam, 1996.


\bibitem[Ye]{Ye}
A. Yekutieli, Dualizing complexes over noncommutative graded
algebras, J. Algebra {\bf 153} (1992), 41-84.

\bibitem[Za]{Za}
A. Zaks, Injective dimension of semiprimary rings, J. Algebra {\bf
13} (1969), 73-86.

\bibitem[Z]{Z}
H. Zassenhaus, The representations of Lie algebras of prime
characteristic, Proc. Glasgow Math. Assoc. {\bf 2} (1954), 1–36.

\end{thebibliography}
\end{document}